\documentclass[10pt]{article}

\usepackage[utf8]{inputenc}
\usepackage[english]{babel}

\usepackage{amssymb,amsmath,amsthm}

\usepackage{lmodern}
\usepackage[T1]{fontenc}

\usepackage{microtype}

\usepackage{comment}

\usepackage[top=1.5cm,bottom=1.5cm,left=1.5cm,right=1.5cm]{geometry}

\usepackage[pdfencoding=unicode, psdextra, citecolor=blue, urlcolor=blue,
linkcolor=blue, colorlinks=true, bookmarksopen=true]{hyperref}

\theoremstyle{plain}
\newtheorem{theorem}{Theorem}

\newtheorem{corollary}{Corollary}
\newtheorem{proposition}{Proposition}
\newtheorem{lemma}{Lemma}
\theoremstyle{definition}

\newtheorem{example}{Example}
\newtheorem{remark}{Remark}

\DeclareMathOperator{\sign}{sign}
\DeclareMathOperator{\supp}{supp}
\DeclareMathOperator*{\argmax}{arg\,max}

\newcommand{\R}{\mathbb{R}}
\newcommand{\Z}{\mathbb{Z}}
\newcommand{\N}{\mathbb{N}}
\newcommand{\SI}{\mathbb{S}}
\newcommand{\T}{{\mathcal T}}
\newcommand{\X}{\mathfrak{X}}
\newcommand{\D}{\,\mathrm{d}}
\newcommand{\rca}{\mathtt{rca}}
\newcommand{\eqdef}{\stackrel{\mathrm{def}}{=}}

\DeclareMathOperator{\E}{\mathcal{E}}

\title{Asymptotic control theory for a closed string II}

\author{Lev Lokutsievskiy \thanks{\textit{lion.lokut@gmail.com}; Steklov Mathematical Institute of Russian Academy of Sciences Moscow, Russia, 119991, Gubkina str.,  8.} \and Alexander Ovseevich \thanks{\textit{ovseev@ipmnet.ru}; Institute for Problems in Mechanics, Russian Academy of Sciences 119526, Vernadsky av., 101/1, Moscow, Russia. }}

\date{}

\begin{document}

\maketitle

\begin{abstract}
We develop an asymptotical control theory for one of the simplest distributed (infinite dimensional) oscillating systems, namely, for a closed string under a bounded load applied to a single distinguished point. We find exact classes of string states that admit complete damping and an asymptotically exact value of the required time. By using approximate reachable sets instead of exact ones, we design a feedback control, which turns out to be asymptotically optimal. The main results are an exact algebraic formula for asymptotic shape of the reachable sets, asymptotically optimal time of motion, and an asymptotically optimal control thus constructed.

\medskip\noindent
\textsc{Keywords}: maximum principle, reachable sets, linear system, string control

\medskip\noindent
\textsc{MSC 2010:} 93B03, 93B07, 93B52.
\end{abstract}

\section{Introduction}
We study a distributed control system governed by
\begin{equation}
\label{basic_system}
  h_{tt}- h_{xx}=u(t)\delta(x),\quad |u|\leq1
\end{equation}
where $h=h(x,t)$ is even and $2\pi$-periodic function on $x$, $u(t)$ is a control and $\delta(x)$ is the Dirac $\delta$-function. Our goal is the design of the minimum time optimal control steering the system from 0 to a given state:
\[
	\begin{gathered}
	T\to\inf\quad\\
	(h,h_t)|_{t=0}=(0,0);
	\qquad
	(h,h_t)|_{t=T}=(h_1,g_1).
	\end{gathered}
\]

The most standard treatment of this system is to regard it as a first order system in the phase space consisting of pairs~$(h,h_t)$ of even and $2\pi$-periodic distributions. Our exposition, however, is based primarily on another form of the system. Let the phase space consist of pairs $(F,f)$, where $F$ is a $2\pi$-periodic distribution (not necessarily even), and $f$ is a real number. The relation with $h$ in the wave equation \eqref{basic_system} is as follows:
\[
	F=h_{x}-h_{t};
	\qquad
	f=\frac1{2\pi}\int_0^{2\pi} h\D x.
\]
Function $h$ can be easily recovered from the data $(F,f)$. Indeed, $h_t$ is even and $h_x$ is odd, hence $h_x(x,t)=\frac12(F(x,t)-F(-x,t))$ and $h(x,t) = f(t) + \int_0^x h_x(y,t)\D y$.

Since $\int_0^{2\pi}h_x\D x=0$, system~\eqref{basic_system} takes the following form
\begin{equation}
\label{eq:basic_system2}
\begin{cases}
	\dot F= -F_x+u\delta\\
	\dot f= -\frac1{2\pi}\int_0^{2\pi} F\D x.
\end{cases}
\end{equation}
If we are given some control $u(t)$, then there exists a unique solution $(F,f)$ to system~\eqref{eq:basic_system2} with given initial conditions $(F,f)|_{t=0}=(0,0)$.

Optimal control problem takes the form
\begin{equation}
\label{problem:main_F_f}
	\begin{gathered}
	T\to\min\qquad\\
	\begin{aligned}
		F|_{t=0}=0;\quad&\quad f|_{t=0}=0;\\
		F|_{t=T}\eqdef F_1=h_{1x}-g_1;\quad&\quad f|_{t=T}\eqdef f_1=\frac1{2\pi}\int_0^{2\pi}h_1\D x.
	\end{aligned}
	\end{gathered}
\end{equation}
We denote by $T(F_1,f_1)$ its value.

\subsection{General strategy}

The basic geometric picture attached to any minimum time control finite dimensional problem is the picture of wave fronts, moving in time $T$, and formed by boundaries $\partial D(T)$ of the sets $D(T)$ reachable from zero in time $T$. The optimal minimum time motion runs in the direction of ``steepest descent'' orthogonal to the wave fronts. This picture is well known as a geometric counterpart of the Bellman ``Dynamic Programming'' based on the Bellman equation. The difficult part in implementation of the ``steepest descent'' strategy is the understanding of the structure of the ``wave fronts'', i.e. the explicit form of the  reachable  sets $D(T)$ as functions on time $T$. In the infinite dimensional case, for the most of end points $(F_1,f_1)$, we have $(F_1,f_1)\not\in\mathrm{cl}\,D(T)$ for $T<T(F_1,f_1)$ and $(F_1,f_1)\in\mathrm{int}\,D(T)$ for $T\ge T(F_1,f_1)$ (see Example~\ref{ex:strange}). So behavior of wave fronts in the infinite dimensional problem~\ref{problem:main_F_f} (despite the problem being linear-convex) is much more complicated than general linear-convex finite dimensional problems.

Surprisingly, this infinite dimensional obstacle can be overcome by the following asymptotic approach, which works much better, than the exact approach described above. In the first place, this approach was developed for the finite dimensional problems and aimed at completely different purposes. Let us describe the approach. It turns out that in a rather general linear control finite dimensional problem one can get an ``explicit'' expression for the reachable  set $D(T)$ in form of an explicit integral formula for the support function of the set. Since $D(T)$ is convex and closed, this gives, in principle, a complete description of the set. Moreover, these sets are well behaved at infinity: There are scaling matrix factors $C(T)$ such that the scaled reachable sets $C(T)D(T)$ tend to a limit set $D_\infty$ as $T\to+\infty.$ This significantly simplifies the structure of the front set at infinity: the front $\partial D(T)$ is approximately $C(T)^{-1}\partial D_\infty$, so that the dependence on time is solely in the scaling matrix $C(T)$ (cf.~\cite{Ovseevich2007}).

This approach provides a solid ground for the ``asymptotical'' control theory of linear systems, where the goal is  the approximate construction of the optimal control ``at infinity''. In the infinite dimensional case, the approach simplifies structure of the boundaries $\partial D(T)$ drastically, since behavior of $C(T)^{-1}\partial D_\infty$ as $T\to\infty$ is very nice and relatively simple.

The remaining hard problem is to extract the structure of the minimum-time control from the explicit knowledge of the support functions of the relevant reachable sets. Morally, this is equivalent to the explicit determination of the corresponding Minkowski functions. In reality this determination is just the first step in construction of the approximate minimum-time control. However, even this preliminary step is quite nontrivial for problem~\eqref{problem:main_F_f}. After this step, we can  define a control which gives the steepest descent with respect to this asymptotic Bellman function. This control appears to be asymptotically optimal.

The aim of this paper is to demonstrate the implementation of this strategy in case of the minimum time control of the closed string. It should be stressed that we do not know the precise scope of success for the above general strategy. It is not clear whether we encounter a good luck, or similar results can be reached in a much greater generality. Anyway,  similar arguments were successfully employed in \cite{Ovseevich2016,Ovseevich2018} for control of a set of linear oscillators, and a simplified damping problem for the ``closed string''.

\subsection{Outline}

There are two main results in this paper. First, Theorem~\ref{thm:T_estimation} gives a precise estimate of the minimal reach time $T(F_1,f_1)$ as $(F_1,f_1)\to\infty$. Second, in Section~\ref{sec:synthesis} we construct asymptotically optimal synthesis for the problem.

The structure of the rest of the paper is roughly as follows.
\begin{enumerate}
	
	\item In Section~\ref{sec:support_function} we get an explicit expression for the support function of the reachable sets $D(T)$.
	
	\item We introduce scaling operators $C(T)$ in Section~\ref{sec:s_d_T_limit} and derive an explicit expression for the limit support function of~$D(T)$.
	
	\item The limit set $D_\infty$ is defined in Section~\ref{sec:d_infty_defn}. We also investigate $D_\infty$ basic properties there.
	
	\item Despite the fact, that we are able to prove a very accurate asymptotic estimate on the optimal time $T(F_1,f_1)$, we first need to derive a very rough factor type estimate on $T(F_1,f_1)$, which is base on $D_\infty$ basic properties (see Section~\ref{sec:T_first_estimate}). Later on we will use the rough estimate to produce the accurate one.

	\item In Section \ref{sec:Minkowski_function} we get an explicit expression for the Minkowski function of the limit reachable set $D_\infty$. This, in turn, requires solution of an auxiliary convex analysis problem (see (\ref{basic_problem})), which is achieved by using Gelfand's representation theory and classical paper of Yosida and Hewitt~\cite{Yosida_Hewitt} on finitely additive measures in section \ref{sec:Gelfand}.
	
	\item On the ground of explicit Minkowski function of the limit reachable set $D_\infty$, we obtain in Section \ref{sec:reach_time} an explicit asymptotic expression for the optimal time $T(F_1,f_1)$.
	
	\item Based on the asymptotic expression for the time of motion we get in Section \ref{sec:synthesis} an explicit admissible control, which allows to bring  the string from complete rest to the given atate in this approximate time.
	
\end{enumerate}

\section{First properties of the reachable sets}

It was shown in~\cite{Ovseevich2018} that reachable states (in any time) of system~\eqref{problem:main_F_f} belongs to $L_\infty(\R/2\pi\Z)\times\R$ even if $F$ is considered as a distribution. So from now on we assume that $F\in L_\infty(\SI)$, where $\SI$ is the 1-dim torus (circle) $S^1=\R/2\pi\Z$.

Let $D(T)$ denote the reachable set of~\eqref{eq:basic_system2} from $(F,f)|_{t=0}=(0,0)$ in time $T$. So, we are interested in time needed to reach a certain point $(F_1,f_1)$, i.e.
\[
	T(F_1,f_1) = \inf\{T\ge0:(F_1,f_1)\in D(T)\}
\]
This infimum is very unlikely to be found explicitly due to complicated structure of $D(T)$. Nonetheless, $D(T)$ has a very nice limiting behavior as $T\to\infty$. In this paper, using this limiting behavior, we find $T$ approximately (see Theorem~\ref{thm:T_estimation}), and the solution allows us to construct asymptotically optimal control.

But before we proceed, let us also underline a strange phenomenon, which appears in this problem for a lot of pairs $(F_1,f_1)$. On one hand, we will prove (see Corollary~\ref{cor:existance_optimal_control}) that for any pair $(F_1,f_1)$ there exists an optimal solution to problem~\eqref{problem:main_F_f}. On the other hand, it usually happens that $(F_1,f_1)\in\mathrm{int}\,D(T)$ for $T=T(F_1,f_1)$ despite the fact that $(F_1,f_1)\not\in D(T)$ for any $T<T(F_1,f_1)$.

\begin{example}
\label{ex:strange}

	Behavior of $F(t)$ in~\eqref{eq:basic_system2} is relatively simple during a period: for $x\in[0;2\pi]$, we have $F(x,2\pi k)=F(x,2\pi(k-1))+u(2\pi k-x)$ for $k\in\N$. Hence $T(F_1,f_1)\ge 2\pi N$ where $N$ is the greatest positive integer that is less than $\|F\|_{L_\infty(\SI)}$. Nonetheless, the behavior of the pair $(F,f)$ in time is much more complicated.
	
	Consider the following function $F_1$: $F_1(x)=1/2$ for $x\in[0;\pi]$ and $F_1(x)=3/2$ for $x\in[\pi;2\pi]$. Let us forget about $f$ for a moment. Then the reach time of $F_1$ from $0$ by first equation in~\eqref{eq:basic_system2} is obviously $3\pi$, since $|u|\le 1$. It is easy to see that $F_1(x)=u(2\pi-x)$ for $x\in[0;\pi]$ and $F(x)=u(2\pi-x)+u(4\pi-x)$ for $x\in [\pi;2\pi]$. Hence $F_1$ can be reached by any control $u(t)\in L_\infty(0;3\pi)$ such that $u(t)=1/2$ for $t\in[\pi;2\pi]$ and $u(t)+u(t+2\pi)=3/2$ for $t\in[0;\pi]$. Therefore there are a variety of admissible controls reaching $F_1$. Simple computation shows that they lead to any $f_1\in[-\frac{5\pi}4,-\frac{7\pi}4]$. Let e.g.\ $f_1=-3\pi/2$. Then $(F_1,f_1)\not\in D(T)$ for any $T<3\pi$, and $(F_1,f_1)\in\mathrm{int}\,D(3\pi)$ since $1/2<1$, $3/2<2$, and $-\frac{5\pi}4<-\frac{3\pi}2<-\frac{7\pi}4$.
	
\end{example}

This example can be generalized by taking any function $F_1$ with $1<\|F_1\|_{L_\infty(\SI)}\not\in\N$. Moreover, if $T\ne\tilde T$, then the Hausdorff distance between $D(T)$ and $D(\tilde T)$ is greater or equal to $1$. This strange phenomenon reflect infinite dimensionality of the state space.

This suggests that  each geometrically intuitive fact from finite dimensional control theory should be carefully reconsidered for the infinite dimension case.

\medskip

We start from deriving first simple properties of $D(T)$

\begin{proposition}
\label{prop:D_T_simple_topology}
	
	For any $T\ge 0$, set $D(T)\subset L_\infty(\SI)\times\R$ is a bounded closed convex set and $D(T)=-D(T)$. If $T\ge T'$, then $D(T)\supset D(T')$.
	
\end{proposition}

\begin{proof}
	
	Let $\E:L_\infty(0;T)\to L_\infty(\SI)\times\R$ denote the end-point map that sends an arbitrary control $u\in L_\infty(0;T)$ to the value of the solution to system~\eqref{eq:basic_system2} with $(F,f)|_{t=0}=(0;0)$ at the final instant $t=T$, $\E(u)=(F,f)|_{t=T}$. It is clear that $D(T)=\E(\mathcal{U})$ where $\mathcal{U}=\{\|u\|_{L_\infty(0;T)}\le 1\}$.
	
	The map $\E$ is linear, and thus convexity and symmetry of the unit ball $\mathcal{U}$ implies convexity and symmetry of $D(T)$.	The map $\E$ is continuous wrt strong topologies on $L_\infty(0;T)$ and $L_\infty(\SI)\times\R$, and thus boundness of the unit ball $\mathcal{U}$ implies boundness of $D(T)$. The map $\E$ is continuous wrt weak$^*$-topologies on $L_\infty(0;T)=L_1^*(0;T)$ and $L_\infty(\SI)\times\R=L_1^*(\SI)\times\R$, and thus weak$^*$-compactness of $\mathcal{U}$ (which follows from Banach--Alaouglu theorem) implies weak$^*$-compactness of $D(T)$. In particular, $D(T)$ is weak$^*$-closed, and so it is strongly closed.
	
	The last property trivially follows from the fact that for any admissible control $u\in L_\infty(0;T')$, we can construct an admissible control in $L_\infty(0;T)$, e.g.
	\[
		\tilde u(t) = \begin{cases}
			0, &t<T-T',\\
			u(t-(T-\tilde T)),& t>T-T',
		\end{cases}
	\]
	that leads to the same endpoint.
\end{proof}

\section{Support function of the reachable sets}\label{sec:support_function}

Sets $D(T)$ are closed and convex, so it is very natural to implement convex analysis. First of all, $(F_1,f_1)\in D(T)$ if and only if $\mu_{D(T)}(F_1,f_1)\le 1$, since $0\in D(T)$ (here $\mu_{D(T)}$ denotes the Minkowski functional). Thus\footnote{Later, we will show (see Corollary~\ref{cor:existance_optimal_control}) that a solution to problem~\eqref{problem:main_F_f} exists, and hence we can write $\min$ instead of $\inf$ here.}
\[
	T(F_1,f_1) = \inf\{T\ge 0:\mu_{D(T)}(F,f)\le 1\}
\]

So our global aim is to find $\mu_{D(T)}$, which is a hard task. Nonetheless, in this section we compute the support function of $D(T)$, which is much simpler.

Introduce the following pairing for $F_1\in L_\infty(\SI)$, $\Phi_1\in L_\infty^*(\SI)$ and $f_1,\phi_1\in \R$:
\[
	\big\langle (F_1,f_1),(\Phi_1,\phi_1) \big\rangle = \langle F_1,\Phi_1 \rangle + 2\pi\,f_1\phi_1,
\]
where $\langle F_1,\Phi_1\rangle$ denotes the pairing between $L_\infty(\SI)$ and $L_\infty^*(\SI)$. For instance, $\|(\Phi_1,\phi_1)\|_{L_1(\SI)\times\R}=\|\Phi_1\|_{L_1(\SI)}+2\pi|\phi|$, since we put $\|(F,f)\|_{L_\infty(\SI)\times\R}=\max\{\|F\|_{L_\infty(\SI)},|f|\}$.

Since $L_1(\SI)$ is weak$^*$-dense in $L^*_\infty(\SI)$ by the Goldstine theorem, it is natural to consider the case $\Phi_1\in L_1(\SI)$. Indeed, if a point $F_1$ in $L_\infty(\SI)$ does not belong to a set $D(T)$ (which is a bounded closed convex set), then it can be strictly separated from $D(T)$ by a linear functional in $L_1(\SI)$ (since $D(T)$ and ball $\{F\in L_\infty(\SI)\ |\ \|F-F_1\|_{L_\infty(\SI)}\le \varepsilon\}$ are both weak$^*$-compact). In other words, values of the support function of $D(T)$ on elements in $L_1(\SI)$ determine $D(T)$ completely.

\begin{proposition}
\label{prop:DT_support_func}
	For any $\Phi_1\in L_1(\SI)$ and $\phi_1\in \R$ we have (recall that $\Phi_1$ is a periodic function)
	\begin{equation}
	\label{eq:s_DT}
	s_{D(T)}(\Phi_1,\phi_1) = \int_0^T|\Phi_1(s)-\phi_1 s|\D s.
	\end{equation}
	If $(\Phi_1,\phi_1)$ defines a supporting half-space at $(F_1,f_1)\in \partial D(T)$, then an optimal control in~\eqref{problem:main_F_f} is given by
	\[
		u(t)\in\sign (\Phi_1(s)-\phi_1 s)\quad\text{where}\quad s=T-t.
	\]
\end{proposition}

\begin{proof}
	If $\Phi_1\in L_1(\SI)$, then
	\[
	\langle F_1,\Phi_1 \rangle=\int_0^{2\pi} F_1(x)\Phi_1(x)\D x.
	\]
	
	To compute the support function of the reachable sets, we consider the adjoint system to the linear part of~\eqref{eq:basic_system2}. Let $\Phi=\Phi(t)$ and $\phi=\phi(t)$. Then the adjoint system takes the form
	\begin{equation}
	\label{basic_system_adjoint}
	\begin{cases}
		\dot{\Phi}=-\Phi_x+\phi\\
		\dot{\phi}=0.
	\end{cases}
	\end{equation}
	The corresponding adjoint Cauchy problem $\Phi(T)=\Phi_1$, $\phi(T)=\phi_1$ has a unique solution:
	\begin{equation}
	\label{explicit}
		\phi(t)\equiv\phi_1
		\qquad\text{and}\qquad
		\Phi(x,t)=\Phi_1(x+T-t)-(T-t)\phi_1.
	\end{equation}
	
	Let $(F(t),f(t))$ and $(\Phi(t),\phi(t))$ be any solutions to systems~\eqref{eq:basic_system2} and~\eqref{basic_system_adjoint}, correspondingly. Then for any $t$ we have
	\[
		\frac{d}{dt}\Big(\big\langle F(t),\Phi(t)\big\rangle + 2\pi\, f(t)\phi(t)\Big) =
		u(t)\Phi(0,t),
	\]
	So we can compute the support fuction $s_{D(T)}$ of the reachable set $D(T)$ from zero in time $T$. Using~\eqref{explicit}, we obtain
	\[
	\begin{aligned}
		s_{D(T)}(\Phi_1,\phi_1)=& \sup_{|u|\le 1} \big(\langle  F(T),\Phi(T)\rangle + 2\pi\, f(T)\phi(T)\big) =
		\sup_{|u|\le 1}\int_0^T \frac{d}{dt}\big(\langle F(t),\Phi(t)\rangle + 2\pi\,f(t)\phi(t)\big)\D t =\\
		=&\sup_{|u|\le 1}\int_0^T u(t) \Phi(0,t)\D t =
		\int_0^T |\Phi(0,t)|dt =
		\int_0^T \left|\Phi_1(T-t)-(T-t)\phi_1\right|dt,
	\end{aligned}
	\]
	where obviously $u=\sign\Phi(0,t) = \sign(\Phi_1(T-t)-(T-t)\phi_1)$.
\end{proof}

\underline{\it From now on we usually omit index $1$} for notations $F_1$, $\Phi_1$, $f_1$, and $\phi_1$ (which have been used to emphasize time independence) and simply write $F\in L_\infty(\SI)$, $\Phi\in L_\infty^*(\SI)$, $f\in\R$, and $\phi\in\R$ for time independent variable when it can not be misinterpreted.

\medskip

Proposition~\ref{prop:DT_support_func} shows that a pair $(\Phi,\phi)$ that defines a supporting half-space to the reachable set $D(T)$ at $(F,f)\in\partial D(T)$ can be used to immediately find an optimal control. So we are also interested in finding the corresponding supporting half-space if $(F,f)\in\partial D(T)$.

\section{Limit behavior of the reachable sets \texorpdfstring{$D(T)$}{D(T)}}
\label{sec:s_d_T_limit}

As it was mentioned, we are interested in the behavior of reachable sets $D(T)$ as $T\to\infty$. It appears that they grow in a very specific way. Introduce scaling operators
\begin{equation}\label{scale}
	C(T)(F,f)= \left(\frac{1}{T}F,\frac{1}{T^2}f\right),
	\qquad
	C(T)^*(\Phi,\phi)= \left(\frac{1}{T}\Phi,\frac{1}{T^2}\phi\right).
\end{equation}
The basic reason for use of these operators is that the scaled reachable sets $C(T)D(T)$ have a nice limit behavior as $T\to\infty$ (it was shown in \cite{Ovseevich2007,Ovseevich2018}). Here we provide an exact estimate of the rate of convergence of the support function $s_{C(T)D(T)}$ as $T\to\infty$:

\begin{proposition}
\label{prop:s_CD_s_D_infty_estimate}

	There exists $c_0>0$ such that for any $\Phi\in L_1(\SI)$, $\phi\in\R$ and $T\ge 2\pi$ we have
	\begin{equation}
	\label{eq:estimate_s_CD_int}
		\left|
			s_{C(T)D(T)}(\Phi,\phi)
			-
			\frac1{2\pi} \int_0^1 \int_0^{2\pi}|\Phi(x)-\phi\tau|\D x\D \tau
		\right| \le
		\frac{c_0}{T} \|(\Phi,\phi)\|_{L_1(\SI)\times\R}
	\end{equation}
	
\end{proposition}

\begin{proof}

	Let $N=[T/2\pi]$. Then for any $\Phi\in L_1(\SI)$ and $\phi\in\R$, we have
	\[
		\begin{aligned}
			s_{C(T)D(T)}(\Phi,\phi) &= s_{D(T)}(C^*(T)(\Phi,\phi)) =
			\frac1T \int_0^T\big|\Phi(s)-\phi\frac sT\big|\D t =\\
			&=\frac1T
			\sum_{k=0}^{N-1}\int_0^{2\pi} \Big|
				\Phi(x)-\phi\frac kN - \phi\frac{x}{2\pi N}
			\Big|\D x
			+\frac1T\int_{2\pi N}^T\big|\Phi(s)-\phi\frac sT\big|\D t.
		\end{aligned}
	\]
	The last term is easy to estimate:
	\[
		\frac1T\int_{2\pi N}^T\big|\Phi(s)-\phi\frac sT\big|\D t \le \frac1T (\|\Phi\|_1 + 2\pi|\phi|).
	\]
	We also prefer to throw away the term $\phi\frac{x}{2\pi N}$ in the sum, which is obviously small as $x\in[0;2\pi]$: $|\phi x / (2\pi N)|\le |\phi|/N$. Hence
	\begin{equation}
	\label{eq:s_CD_intermediate_int}
		\left|
			s_{C(T)D(T)}(\Phi,\phi) -
			\frac1T
			\sum_{k=0}^{N-1}\int_0^{2\pi} \Big|
				\Phi(x)-\phi\frac kN
			\Big|\D x
		\right| \le \frac1T (\|\Phi\|_1 + 4\pi|\phi|)
	\end{equation}
	
	Now for the integral in~\eqref{eq:estimate_s_CD_int}, we consider the Riemannian sum constructed by the division of $[0;1]$ into $N$ equal parts:
	\[
		\mathcal{R}_N(\xi)\eqdef\frac1{2\pi N}\sum_{k=0}^{N-1}\int_0^{2\pi} \Big|
			\Phi(x)-\phi\xi_k
		\Big|\D x
	\]
	where $\xi=(\xi_0,\ldots,\xi_{N-1})$ and $\xi_k\in[\frac kN;\frac{k+1}N]$. For $\xi_k=k/N$, we obtain the sum in~\eqref{eq:s_CD_intermediate_int}. Therefore, since both the integral and $\mathcal{R}_N(k/N)$ belong to the interval $[\inf_{\xi}\mathcal{R}_N(\xi),\sup_{\xi}\mathcal{R}_N(\xi)]$, we obtain
	\[
		\left|
			\frac1{2\pi} \int_0^1 \int_0^{2\pi}|\Phi(x)-\phi\tau|\D x\D \tau
			-
			\frac1{2\pi N}
			\sum_{k=0}^{N-1}\int_0^{2\pi}
			\Big|
				\Phi(x)-\phi\frac kN
			\Big|\D x
		\right|
		\le
		\frac{|\phi|}{N} \le \frac{4\pi|\phi|}{T}
	\]
	Finally, let us replace factor at the sum $\frac1{2\pi N}$ by $\frac1T$. Since $0\le \frac{1}{2\pi N}-\frac1T\le \frac1{NT}$ and
	\[
		\sum_{k=0}^{N-1}\int_0^{2\pi}
		\Big|
			\Phi(x)-\phi\frac kN
		\Big|\D x
		\le
		N(\|\Phi\|_{L_1(\SI)}+2\pi|\phi|),
	\]
	this change can increase the error by $\frac1T(\|\Phi\|_{L_1(\SI)}+2\pi|\phi|)$ at most. Thus $c_0=5$ fits the proposition statement.
\end{proof}

\section{Limit reachable set \texorpdfstring{$D_\infty$}{D∞}}
\label{sec:d_infty_defn}

Proposition~\ref{prop:s_CD_s_D_infty_estimate} shows that $s_{C(T)D(T)}$ has the following uniform limit on any bounded set in $L_1(\SI)\times\R$ as $T\to\infty$:
\[
	s_{C(T)D(T)}\rightrightarrows \frac1{2\pi} \int_0^1 \int_0^{2\pi}|\Phi(x)-\phi\tau|\D x\D \tau.
\]
Thus, the rhs is a convex positively homogeneous function on $L_1(\R)\times\R$, and hence it can be considered as a support function of the limiting (in some sense) set
\[
	D_\infty=\lim_{T\to\infty}C(T)D(T),
\]
where
\begin{equation}
\label{eq:D_infty_defn}
	D_\infty \eqdef \Big\{
		(F,f)\in L_\infty(\SI)\times\R\quad\Big|\quad\forall (\Phi,\phi)\in L_1(\SI)\times\R\quad\text{we have}\quad \langle (F,f),(\Phi,\phi)\rangle\le \frac1{2\pi} \int_0^1 \int_0^{2\pi}|\Phi(x)-\phi\tau|\D x\D \tau
	\Big\},
\end{equation}
Note that we have defined $D_\infty$ using only $\Phi\in L_1(\SI)$, and points $\Phi\in L^*_\infty(\SI)\setminus L_1(\SI)$ remain unused.

Now, we find the value of the support function of the limit set $D_\infty$  at $\Phi\in L_1(\SI)$. We claim that
\begin{equation}
\label{eq:s_D_infty}
	\forall (\Phi,\phi)\in L_1(\SI)\times\R\qquad
	s_{D_\infty}(\Phi,\phi)=\frac1{2\pi} \int_0^1 \int_0^{2\pi}|\Phi(x)-\phi\tau|\D x\D \tau.
\end{equation}
Indeed, let us denote by $\tilde s(\Phi,\phi)$ the rhs of~\eqref{eq:s_D_infty} and then prove that it coincides with the support function of $D_\infty$. Function $\tilde s$ is convex and continuous wrt strong topology on $L_1(\SI)$. Hence $\tilde s$ is closed under weak topology by the Banach separation theorem. The spaces $L_1(\SI)\times \R$ and $L_\infty(\SI)\times\R$ with weak- and weak$^*$- topologies correspondingly form a pair of spaces in duality (see~\cite{Tikhomirov_Magaril}). Thus $\tilde s=\tilde s^{**}$ by the Fenchel--Moreau theorem\footnote{Here $\tilde s^*$ and $\tilde s^{**}$ are conjugate and double Legendre--Young--Fenchel conjugate functions defined on $L_\infty(\SI)\times\R$ and $L_1(\SI)\times \R$, correspondingly.}. It remains to note that $\tilde s^*$ is the indicator function of $D_\infty$ by definition of $D_\infty$.

So we are able to easily find values of the support function $s_{D_\infty}$ at $\Phi\in L_1(\SI)$. Unfortunately, values of $s_{D_\infty}$ for $\Phi\in L^*_\infty(\SI)$ are not so easy to compute, and we postpone the computation until section~\ref{sec:Gelfand}, where these values are found by using the Gelfand transform.

\begin{proposition}
\label{prop:s_D_infty_norm}
	
	Function $s_{D_\infty}$ is a norm on $L_1(\SI)\times\R$, which is equivalent to $\|\cdot\|_{L_1(\SI)\times\R}$.
	
\end{proposition}

\begin{proof}
	
	Obviously, $s_{D_\infty}$ is a convex, non-negative, positively homogeneous function, and $s_{D_\infty}(-\Phi,-\phi)=s_{D_\infty}(\Phi,\phi)$. Hence it defines a norm on $L_1(\SI)\times\R$. We compare it with the standard one $\|\cdot\|_{L_1(\SI)\times\R}$.
	
	First, for any $\Phi\in L_1(\SI)$ and $\phi\in\R$, we have
	\[
		s_{D_\infty}(\Phi,\phi) \le |\phi| + \frac{1}{2\pi}\|\Phi\|_{L_1(\SI)}.
	\]
	
	We show now, that there is a constant $c>0$ such that $s_{D_\infty}(\Phi,\phi)\ge c(\|\Phi\|_{L_1(T)}+|\phi|)$. If $\phi=0$, then $2\pi s_{D_\infty}(\Phi,0)=\|\Phi\|_{L_1(\SI)}$. Suppose $\phi\ne 0$. Since $s_{D_\infty}$ is even positively homogeneous functional, we can put $\phi=1$ without loss of generality. Using Fubini's theorem we may change  the order of integration in~\eqref{eq:s_D_infty}. There are three cases:
	\[
		\begin{aligned}
			\Phi(x)\le 0 \quad\Rightarrow\quad& \int_0^1|\Phi(x)-\tau|\D \tau = |\Phi(x)| + \frac12;\\
			\Phi(x)\ge 1 \quad\Rightarrow\quad& \int_0^1|\Phi(x)-\tau|\D \tau = \Phi(x)-\frac12 \ge \frac14\Phi(x)+\frac14;\\
			\Phi(x)\in[0;1] \quad\Rightarrow\quad& \int_0^1|\Phi(x)-\tau|\D \tau = \Phi^2(x)-\Phi(x)+\frac12 \ge \frac18\Phi(x)+\frac18;
		\end{aligned}
	\]
	Therefore,
	\[
		s_{D_\infty}(\Phi,1)\ge \frac18 (\|\Phi\|_{L_1(\SI)} + 1).
	\]
	
\end{proof}

\begin{corollary}
	\label{cor:D_infty_topology}
	
	Set $D_\infty\subset L_\infty(\SI)\times\R$ is a bounded closed convex set, $D_\infty=-D_\infty$, and $0\in\mathrm{int}\,D_\infty$.
	
\end{corollary}

\section{Estimates of the reach time}
\label{sec:T_first_estimate}

In this section we show that the initial problem~\eqref{problem:main_F_f} has an optimal solution for any end point. This result is closely related to some basic estimations on $T(F,f)$.

\begin{proposition}
\label{prop:D_infty_CTDT_approx}

	There exists $T_1>0$ such that $T\ge T_1$ implies
	\[
		\Big(1-\frac{T_1}T\Big)D_\infty
		\subset
		 C(T)D(T)
		\subset
		\Big(1+\frac{T_1}T\Big)D_\infty.
	\]
	
\end{proposition}

\begin{proof}
	
	Since $s_{D_\infty}$ is equivalent to the standard norm on $L_1(\SI)\times\R$ by Proposition~\ref{prop:s_D_infty_norm}, using Proposition~\ref{prop:s_CD_s_D_infty_estimate} we obtain
	\[
		\left|
			s_{C(T)D(T)}(\Phi,\phi) - s_{D_\infty}(\Phi,\phi)
		\right| \le
		\frac{c_0}T\|(\Phi,\phi)\|_{L_1(\SI)\times\R} \le
		\frac{cc_0}{T} s_{D_\infty}(\Phi,\phi)
	\]
	for some constant $c>0$ and any $(\Phi,\phi)\in L_1(\SI)\times\R$. Consequently,
	\[
		\big(1-\frac{cc_0}{T}\big) s_{D_\infty}(\Phi,\phi)
		\le
		s_{C(T)D(T)}(\Phi,\phi)
		\le
		\big(1+\frac{cc_0}{T}\big) s_{D_\infty}(\Phi,\phi).
	\]
	Put $T_1=cc_0$. Since both $D_\infty$ and $C(T)D(T)$ are closed convex bounded sets, they are weak$^*$ compact. Hence if $T\ge T_1$, then the written inequalities on $s_{D_\infty}$ and $s_{C(T)D(T)}$ are equivalent to the Proposition statement by the bipolar theorem (applied to the pair of spaces in duality, namely $L_1(\SI)\times\R$ and $L_\infty(\SI)\times\R$).

\end{proof}

\begin{corollary}
\label{cor:Minkowski_functional_approx}
	For $T\ge T_1$ and any $(F,f)\in L_\infty(\SI)\times\R$ we have
	\[
		\Big(1-\frac{T_1}{T}\Big)\mu_{C(T)D(T)}(F,f)
		\le
		\mu_{D_\infty}(F,f)
		\le
		\Big(1+\frac{T_1}{T}\Big)\mu_{C(T)D(T)}(F,f)
	\]
	
\end{corollary}

\begin{proof}
	
	This immediately follows from Proposition~\ref{prop:D_infty_CTDT_approx}
	
\end{proof}

\begin{corollary}
\label{cor:existance_optimal_control}

	There exists an optimal solution to the initial problem~\eqref{problem:main_F_f} for any given end point.

\end{corollary}

\begin{proof}
	
	First, let us show, that for any $(F,f)\in L_\infty(\SI)\times\R$ there exists an admissible control that leads to $(F,f)$. In other words, we claim that
	\[
		\bigcup_{T\ge 0}D(T) = L_\infty(\SI)\times\R.
	\]
	Indeed, for $T\ge 2T_1$, set $C(T)D(T)$ contain a ball of some fixed radius $r$ by Corollary~\ref{cor:D_infty_topology} and Proposition~\ref{prop:D_infty_CTDT_approx}. Hence $D(T)$ contains a ball of radius $Tr/(2T_1)$ for $T\ge \max\{1,2T_1\}$.
	
	Second, let us now prove the existence of an optimal control. Fix $(F,f)$. Let $u_n\in L_\infty(T_n)$, $n\in\N$, be a minimizing sequence, i.e.\ $\E(u_n)=(F,f)$ for all $n$ and $T_n\to T(F,f)$ as $n\to\infty$. Chose $\hat T>T_n$ for all $n$, and put $\tilde u_n(t)=0$ for $t\in(0;\hat T-T_n)$ and $\tilde u_n(t)=u_n(t+(\hat T-T))$ for $t\in(\hat T-T_n,\hat T)$. Then $\E(\tilde u_n)=(F,f)$. Obviously, $\tilde u_n\in L_\infty(0;\hat T)$ and $\|\tilde u_n\|\le 1$. Therefore, there exists a subsequence that  weak$^*$-converges to a control $\tilde u\in L_\infty(0;\hat T)$. Since $\E$ is continuous wrt weak$^*$ topologies, we have $\E(\tilde u)=(F,f)$. Moreover $u(t)=0$ for $t<\hat T-T(F,f)$. Hence the control $u(t)=\tilde u(t-(\hat T-T(F,f)))$, $u\in L_\infty(0;T(F,f))$ is optimal.
	
\end{proof}

\begin{corollary}
\label{cor:reach_time_simple_estimation}

	There exists $c_1>0$ such that for all $(F,f)$ with $\|(F,f)\|_{L_\infty(\SI)\times\R}\ge c_1$, minimal reach time $T(F,f)$ satisfies
	\[
		\frac1{c_1}\big(\|F\|_{L_\infty(\SI)} + |f|^{1/2}\big) \le
		T(F,f) \le
		c_1\big(\|F\|_{L_\infty(\SI)} + |f|^{1/2}\big).
	\]

\end{corollary}

\begin{proof}

	Obviously, $\mu_{D(T)}(F,f)=\mu_{C(T)D(T)}(C(T)(F,f))\le 1$ and $\mu_{C(\tilde T)D(T')}(C(T')(F,f))>1$ for any $T'<T$. Function $\mu_{D_\infty}(C(T')(F,f))$ is continuous in $T'$. Hence, if $T>2T_1$, then by Corollary~\ref{cor:Minkowski_functional_approx} we obtain $\frac12 \le	\mu_{D_\infty}(C(T)(F,f)) \le \frac32$.	Since $\mu_{D_\infty}$ is equivalent to standard norm, it follows that
	\[
		1/c \le
		\frac1T\max\Big\{\|F\|_{L_\infty(\SI)},\frac fT\Big\}\le
		c
	\]
	for some constant $c>0$. The right inequality implies both $T\ge \|F\|_{L_\infty(\SI)}/c$ and $T\ge \sqrt{|f|/c}$. The left inequality implies $T\le c\|F\|_{L_\infty(\SI)}$ or $T\le \sqrt{c|f|}$. Summarizing, for $T\ge2T_1$ we obtain
	\[
		\frac1{2c}\big(\|F\|_{L_\infty(\SI)} + |f|^{1/2}\big)\le
		T\le
		c\big(\|F\|_{L_\infty(\SI)} + |f|^{1/2}\big),
	\]
	since $c\ge 1$. Set $D(2T_1)$ is bounded, and hence the inequality $T(F,f)\ge 2T_1$ is satisfied for all $(F,f)$ with sufficiently large norm. It remain to put $c_1=\max\{2c,\mathrm{diam}\, D(2T_1)\}$.
\end{proof}

\section{Gelfand representation}
\label{sec:Gelfand}

Corollaries~\ref{cor:Minkowski_functional_approx} and~\ref{cor:reach_time_simple_estimation} show that for large enough $(F,f)$, the minimal reach time $T(F,f)=\min\{T:\mu_{D(T)}(F,f)\le 1\}$ is very close to a solution of the equation $\mu_{D_\infty}(C(T)(F,f))=1$, and allows us to estimate the error.

This makes important to compute $\mu_{D_\infty}$ explicitly. The only remaining obstacle is that we do know values of the support function $s_{D_\infty}(\Phi,\phi)$ for $\Phi\in L_1(\SI)$ (see~\eqref{eq:s_D_infty}), but we do not know its values for $\Phi\in L_\infty^*(\SI)\setminus L_1(\SI)$.

Moreover, it is necessary to use the dual space $L_\infty^*(\SI)$ for describing $\partial D_\infty$. Indeed, any point at the boundary $\partial D_\infty$ has a supporting half-space, but it is determined by a covector $(\Phi,\phi)\in L_\infty^*(\SI)\times\R$, which does not belong to $L_1(\SI)\times \R$ in general\footnote{Structure of supporting half-spaces here is similar to the structure of supporting half-spaces to the unit ball $\{\|F\|_{L_\infty(\SI)}\le 1, F\in L_\infty(\SI)\}$, which are not always determined by $L_1(\SI)$ functions.}.

The space $L_\infty^*(\SI)$ consists of finitely additive signed measures on $\SI$ that are absolutely continuous wrt Lebesgue measure on $\SI$. Finitely additive measures are not so very well behaved as classical countable additive ones. Fortunately, we may avoid considerations of finitely addivitive measures by using Gelfand's representation of Banach algebras. According to \cite[Theorem 4.3]{Yosida_Hewitt}, the space $L_\infty(\SI)$ is isomorphic (as a Banach space) to the space of continuous functions $C(\X)$ on a topological space~$\X$ of non-trivial characters on $L_\infty(\SI)$, namely
\[
	\X = \Big\{
		\chi\in L_\infty^*(\SI)
		\quad\Big|\quad
		\chi:L_\infty(\SI)\to\R\ \ \text{is a non-trivial algebra homomorphism}
	\Big\}.
\]
The topology on $\X$ is the weak$^*$-topology induced by inclusion $\X\subset L^*_\infty(\SI)$. Obviously, the space $L_\infty^*(\SI)$ is weak$^*$-Hausdorff, and $\X$ is a weak$^*$-closed subset of the unit ball. Hence $\X$ is a compact Hausdorff topological space by the Banach--Alaoglu theorem. The map $L_\infty(\SI)\to C(\X)$ is called Gelfand representation and it is an isomorphism of Banach spaces and even algebras~\cite[Theorem 4.3]{Yosida_Hewitt}. The image of a function $F\in L_\infty(\SI)$ under the isomorphism we denote by $\bar F\in C(\X)$, namely
\[
	\bar F(\chi)\eqdef\langle\chi,F\rangle.
\]

We are interested in the Gelfand representation, since the space $L_\infty^*(\SI)$ is isomorphic to the dual space $C^*(\X)$, which is much nicer than $L^*_\infty(\SI)$ as $\X$ is a compact Hausdorff topological space. Indeed, $C^*(\X)$ coincides with the space $\rca(\X)$ of (finite, countable additive) signed regular Borel measures on $\X$ by the Riesz--Markov--Kakutani representation theorem. For any linear functional $\Phi\in L_\infty^*(\SI)$ (e.g.\ $\Phi\in L_1(\SI)$), we denote by $\bar\Phi\in\rca(\X)$ its image under this isomorphism. Obviously, $\|\bar\Phi\|_{C(\X)}=\|\Phi\|_{L^\infty(\SI)}$.

Regular Borel measures on a (not perfectly nice) topological space $\X$ are much more easy to work with than finitely additive measures on $\SI$. For example, the following proposition allows us to extend formula~\eqref{eq:D_infty_defn} from $L_1(\SI)\times \R$ to $L_\infty^*(\SI)\times\R$, and write $s_{D\infty}$ explicitly for all $\Phi\in L_\infty^*(\SI)$ in terms of Gelfand representation.

\begin{proposition}

	For any $(\Phi,\phi)\in L_\infty^*(\SI)\times\R$, we have
	\begin{equation}
	\label{eq:s_D_infty_X}
		s_{D_\infty}(\Phi,\phi) =
		\frac1{2\pi}\int_0^1\int_\X \big|\bar\Phi - \varphi\tau\bar1\big|\D \tau.
	\end{equation}

\end{proposition}

\begin{proof}
	
	Denote by $\tilde s(\Phi,\phi)$ the rhs of~\eqref{eq:s_D_infty_X}. First, let us prove that $s_{D_\infty}$ and $\tilde s$ define norms on $L_\infty^*(\SI)\times \R$ that are both equivalent to $\|\cdot\|_{L_\infty^*(\SI)\times \R}$. Obviously, $s_{D_\infty}\sim \|\cdot\|_{L_\infty^*(\SI)\times \R}$ by Proposition~\ref{prop:s_D_infty_norm}. Moreover, $\tilde s(\Phi,\phi)\le |\phi| + \frac{1}{2\pi}\|\bar\Phi\|_{C^*(\X)}$. To show the opposite inequality, we use the Lebesgue decomposition theorem, which implies that $\bar\Phi=\bar\Phi^{ac}+\bar\Phi^s$, where measure $\bar\Phi^{ac}$ is absolutely continuous wrt $\bar1$, $\bar\Phi^{ac}\ll\bar 1$, and $\bar\Phi^s$ is singular wrt $\bar 1$, $\bar \Phi^s\perp\bar 1$. Hence $|\bar\Phi - \varphi\tau\bar1| = |\bar\Phi^{ac} - \varphi\tau\bar1| + |\bar\Phi^s|$, and
	\[
		\begin{aligned}
			2\pi s_{D_\infty}(\Phi,\phi) = &
			\|\bar\Phi^s\|_{C(\X^*)} + \int_0^1\int_0^{2\pi} |\Phi^{ac}(x)-\phi\tau|\D x\D \tau = \\
			=&\|\Phi^s\|_{L_\infty^*(\SI)} + 2\pi s_{D_\infty}(\Phi^{ac},\phi)\ge
			\|\Phi^s\|_{L_\infty^*(\SI)} + c(\|\Phi^{ac}\|_{L_1(\SI)} + |\phi|).
		\end{aligned}
	\]
	
	So norms $s_{D_\infty}$ and $\tilde s$ are both equivalent to $\|\cdot\|_{L_\infty^*(\SI)\times\R}$ and coincide on $L_1(\SI)\times \R$. Hence $s_{D_\infty}\equiv \tilde s$ on $L_\infty^*(\SI)\times \R$ by the lemma below.
		
\end{proof}

The following lemma looks pretty standard, but we were unable to find an exact reference and decided to give a proof.

\begin{lemma}
\label{prop:bidual_norm_by_weak_star_limit}

	Let $Y$ be a Banach space and $\zeta_1,\zeta_2:Y^{**}\to\R_+$ are two norms that are both equivalent to $\|\cdot\|^{**}$. If $\zeta_1(y)=\zeta_2(y)$ for all $y\in Y$, then $\zeta_1(w)=\zeta_2(w)$ for all $w\in Y^{**}$.

\end{lemma}

\begin{proof}

	Sets $\{w\in Y^{**}:\zeta_{1,2}(w)\le a\}$ are convex and closed for all $a\in\R$, since functions $\zeta_{1,2}$ are convex and continuous. Moreover, these sets are bounded, since $\|\cdot\|^{**}$ is continuous wrt $\zeta_{1,2}$. Hence, norms $\zeta_{1,2}$ are weak$^*$-closed by Banach-Alaoglu theorem. Therefore $\zeta_1\equiv\zeta_2$ iff $\zeta_1^*\equiv\zeta_2^*$ by the Fenchel–Moreau theorem (applied to the pair of spaces $(Y^*,Y^{**})$ in duality). We have $\zeta_{1,2}^*:Y^*\to \R$, where $\zeta_{1,2}^*=\delta_{B_{1,2}}$ for $B_{1,2}=\{z\in Y^*:\forall w\in Y^{**}\ \langle z,w\rangle\le\zeta_{1,2}(w)\}$.
	
	Here is the trick: consider restriction of the biconjugate functions $\delta^*_{B_{1,2}}=\zeta_{1,2}$ to $Y$. Sets $B_{1,2}\in Y^*$ are bounded and closed. Thus they are both weak$^*$-closed in $Y^*$ by the Banach--Alaouglu theorem. So the functions $\delta_{B_{1,2}}$ are also weak$^*$-closed on $Y^*$, and we can apply the Fenchel–Moreau theorem for the pair of spaces $(Y,Y^*)$ in duality, which implies that $\delta_{B_1}\equiv \delta_{B_2}$, since $\delta^*_{B_1}|_{Y}\equiv \delta^*_{B_2}|_Y$.

\end{proof}

\section{The limit Minkowski function}\label{sec:Minkowski_function}
\begin{theorem}
\label{thm:main}
	The value of the Minkowski function $\mu_{D_\infty}(F,f)$ of the limit reachable set $D_\infty$ is given by
	\[
		\mu_{D_\infty}(F,f) = 2\pi\max\left\{
			\|F\|_{L_\infty(\SI)},\
			\big|\langle F\rangle + 2f\big| + \sqrt{(\langle F\rangle + 2f)^2+\langle F^2\rangle}
		\right\}
	\]
	where $\langle G\rangle=\frac1{2\pi}\langle G,1\rangle$ stands for the mean value of $G\in L_\infty(\SI)$.
	
	If $\mu_{D_\infty}(F,f)>2\pi\|F\|_{L_\infty(\SI)}$, then there exists a unique pair $(\Phi,\phi)\in L_\infty^*(\SI)\times\R$ such that
	\[
		\langle F,\Phi\rangle + 2\pi\,f\phi = \mu_{D_\infty}(F,f)
		\qquad\text{and}\qquad
		s_{D_\infty}(\Phi,\phi)=1
	\]
	Moreover, in this case,
	\begin{equation}
	\label{eq:optimal_phi_L_one_case}
		\phi = \frac{4\mu^2}{\mu^2+4\pi^2\langle F^2 \rangle}\sign\big(\langle F\rangle + 2f\big),
		\qquad\text{and}\qquad
		\Phi(x)=\frac\pi\mu |\phi| F(x) + \frac12\phi
	\end{equation}
	where $\mu=\mu_{D_\infty}(F,f)$ for short.
	
\end{theorem}

\begin{proof}
	
	First consider the following auxiliary problem on $\Phi\in L_\infty^*(\SI)$ and $\phi\in\R$
	\begin{equation}\label{basic_problem}
		\left\{\begin{gathered}
			J(\Phi,\phi)=\big\langle (F,f),(\Phi,\phi)\big\rangle \to \max_{(\Phi,\phi)}\\
			s_{D_\infty}(\Phi,\phi)\le 1.
		\end{gathered}\right.
	\end{equation}
	Denote by $V(F,f)$ the auxiliary problem value.

	\begin{lemma}
	\label{lm:functional_V}

		$\mu_{D_\infty}\equiv V$ and hence functional $V$ is convex, positively-homogeneous, and continuous.

	\end{lemma}

	\begin{proof}

		Put
		\[
			D_\infty^\circ \eqdef
			\Big\{
				(\Phi,\phi)\in L^*_\infty(\SI)\times\R\quad\Big|\quad s_{D_\infty}(\Phi,\phi)\le 1
			\Big\}.
		\]
		Suppose $\mu_{D_\infty}(F,f)=1$ without loss of generality. Obviously, for any $(\Phi,\phi)\in D_\infty^\circ$, we have $\langle (F,f),(\Phi,\phi)\rangle\le 1$ as $(F,f)\in D_\infty$, and hence $V(F,f)\le 1$. Let us show that $V(F,f)\ge1$. The problem is homogeneous, hence $\langle (F,f),(\Phi,\phi)\rangle\le V(F,f)s_{D_\infty}(\Phi,\phi)$ for any $\Phi\in L_\infty^*(\SI)$ and $\phi\in\R$. The $D_\infty$ is a closed convex set by Proposition~\ref{cor:D_infty_topology}, and hence $(F,f)\in V(F,f) D_\infty$ by bipolar theorem. Therefore $1=\mu_{D_\infty}(F,f)\le V(F,f)$.
		
		It remains to note that $\mu_{D_\infty}$ a norm on $L_\infty(\SI)\times\R$, which is equivalent to the standard one by Corollary~\ref{cor:D_infty_topology}.


	\end{proof}

	So, Minkowski functional $\mu_{D_\infty}$ coincides with the auxiliary problem value $V$, which we are able to compute explicitly. We start with the simple

	\begin{lemma}

		$V(F,f)\ge 2\pi\|F\|_{L_\infty(\SI)}$.

	\end{lemma}

	\begin{proof}

		Indeed, if $F=0$ then there is nothing to be proved. Let $F\ne 0$. Then, for any $0<\varepsilon<\|F\|_{L_\infty(\SI)}$, consider the following set
		\[
			E_\varepsilon = \{x\in\SI: |F(x)|>\|F\|_{L_\infty(\SI)}-\varepsilon\}
		\]
		of non-zero measure $|E_\varepsilon|$. Put $\phi=0$ and
		\[
			\Phi_\varepsilon(x) =
			\begin{cases}
				0,&\text{if}\ x\not\in E_\varepsilon;\\
				2\pi \sign F(x)/|E_\varepsilon|,&\text{if}\ x\in E_\varepsilon.
			\end{cases}
		\]
		Then $\Phi_\varepsilon\in L_\infty(\SI)\subset L_1(\SI)$, $s_{D_\infty}(\Phi,0)=1$ and
		\[
			\langle F,\Phi\rangle+2\pi f\phi = \frac{2\pi}{|E_\varepsilon|} \int_{E_\varepsilon}|F(x)|\D x \ge 2\pi(\|F\|_{L_\infty(\SI)}-\varepsilon).
		\]
		Since $\varepsilon$ can be chosen arbitrary close to $0$, we obtain $V(F,f)\ge 2\pi\|F\|_{L_\infty(\SI)}$.

	\end{proof}
	
	Note that, for any $(F,f)$, there exists a solution to the auxiliary problem. Indeed, set $D_\infty^\circ$ is a bounded closed convex set, hence it is weak$^*$-compact in $L_\infty^*(\SI)\times\R$ by Banach--Alaoglu theorem, and the functional $J$ is obviously weak$^*$-continuous. Hence, in order to compute $V(F,f)$ we apply Karush–Kuhn–Tucker (KKT) theorem, which (together with Slater's condition) implies that there exists $\mu\ge 0$ such that the Lagrangian
	\[
		\mathcal{L}(\Phi,\phi) = -\langle (F,f),(\Phi,\phi)\rangle + \mu\, s_{D_\infty}(\Phi,\phi)
	\]
	attains its global minimum at the solution to the auxiliary problem and the complementary slackness condition $\mu(s_{D_\infty}(\Phi,\phi)-1)=0$ is fulfilled. Slightly abusing notations, we denote by $(\Phi,\phi)$ a solution to the auxiliary problem. Then
	\[
		0\in \partial\mathcal{L}
		\qquad\Longleftrightarrow\qquad
		(F,f) \in \mu\, \partial s_{D_\infty}(\Phi,\phi).
	\]
	Since
	\[
		\partial s_{D_\infty}(\Phi,\phi) = \argmax_{(F_1,f_1)\in D_\infty}\langle (F_1,f_1),(\Phi,\phi)\rangle,
	\]
	we obtain $\mu=V(F,f)$.
	
	So we need to prove that (i) if $\mu>2\pi\|F\|_{L_\infty(\SI)}$, then
	\begin{equation}
	\label{eq:mu_second}
		\mu=2\pi\left[\big|\langle F\rangle + 2f\big| + \sqrt{(\langle F\rangle + 2f)^2+\langle F^2\rangle}\right]
	\end{equation}
	and (ii) if rhs of~\eqref{eq:mu_second} if strictly greater than $2\pi\|F\|_{L_\infty(\SI)}$, then $\mu$ equals rhs of~\eqref{eq:mu_second}.
	\medskip
	
	First, we suppose that $\mu>2\pi\|F\|_{L_\infty(\SI)}$. Let us show that in this case, the auxiliary problem has the solution $(\Phi,\phi)$ given in the theorem statement.
	
	Now we find $\partial_\Phi \mathcal{L}$. Denote by $\bar\Psi_+$ and $\bar\Psi_-$ the positive and negative components of a measure $\bar\Psi\in\rca(\X)$, i.e.\ $\bar\Psi=\bar\Psi_+-\bar\Psi_-$, $\bar\Psi_+\ge0$, $\bar\Psi_-\ge 0$, and $\bar\Psi_+\perp\bar\Psi_-$. Then
	\[
		\partial|\bar\Psi| = \Big\{
			\bar F\in C(\X)
			\quad\Big|\quad
			\|\bar F\|_{C(\X)}\le 1,\ \
			\bar F|_{\supp\bar\Psi_+}\equiv 1,\ \ \text{and}\ \ \bar F|_{\supp\bar\Psi_-}\equiv-1
		\Big\}\subset C(\X).
	\]
	Note that subdifferential $\partial|\bar\Psi|\subset C(\X)$ may be empty, namely if $\supp\bar\Psi_+\cap \supp\bar\Psi_-\ne\emptyset$. Nonetheless, solution to auxiliary problem exists, hence KKT conditions dictates that the corresponding subdifferential of $\mathcal L$ is not empty and contains 0. So,
	\[
		\bar F \in \mu\,\partial_\Phi\mathcal{L} =
		\frac\mu{2\pi} \int_0^1 \partial_{\bar\Phi} |\bar\Phi -\phi\tau\bar 1|\D \tau.
	\]
	or
	\[
		\forall \chi\in\X
		\qquad
		\bar F(\chi) \in \frac\mu{2\pi}
		\int_0^1 S(\chi,\tau)\D \tau
	\]
	where
	\[
		S(\chi,\tau) = \begin{cases}
			1,&\text{if }\chi\in \supp (\bar\Phi - \phi\tau\bar1)_+;\\
			-1,&\text{if }\chi\in \supp (\bar\Phi - \phi\tau\bar1)_-;\\
			[-1;1],&\text{otherwise}.
		\end{cases}
	\]
	
	We claim that $\phi\ne0$. Indeed, if $\phi=0$ then for any $\chi\in\supp \bar\Phi$ we have $|\bar F(\chi)|>\|F\|_{L_\infty(\SI)}=\|\bar F\|_{C(\X)}$ as $\mu>2\pi\|F\|_{L_\infty(\SI)}$, which is impossible. Thus $\supp\bar\Phi=\emptyset$ and the pair $(\Phi=0,\phi=0)$ does not satisfy the complementary slackness condition. Hence, $\phi\ne 0$. Moreover, we claim that the measure $\bar\Phi$ must be absolutely continuous wrt $\bar1$, $\bar\Phi\ll\bar 1$. Indeed, using the Lebesgue decomposition theorem we can write $\bar\Phi = \bar\Phi^{ac}+\bar\Phi^s$, where $\bar\Phi^{ac}\ll\bar 1$ and $\bar\Phi^s\perp\bar1$. For all $\tau$ we have $(\bar\Phi-\phi\tau\bar1)_\pm = \bar\Phi^s_\pm + (\bar\Phi^{ac}-\phi\tau\bar1)_\pm$, and hence $\supp(\bar\Phi-\phi\tau\bar1)_\pm\supset \supp \bar\Phi^s_\pm$. Thus for all $\chi\in\supp\bar\Phi^s$ we have $|\bar F(\chi)|>\|\bar F\|_{C(\X)}$, which is impossible. So $\supp\bar\Phi^s=\emptyset$ and the measure $\bar\Phi$ belongs to $L_1(\X,\bar 1)$. Therefore $\bar\Phi$ can be considered as a measurable function on $\X$, and for $\bar1$-a.e.\ $\chi$ we have
	\[
		S(\chi,\tau) = \begin{cases}
			1,&\text{if}\quad\bar\Phi(\chi) > \phi\tau;\\
			-1,&\text{if}\quad\bar\Phi(\chi)< \phi\tau;\\
			[-1;1],&\text{if}\quad\bar\Phi(\chi) = \phi\tau;.
		\end{cases}
	\]
	
	Put $\bar\tau_*(\chi) = \bar\Phi(\chi)/\phi$ for $\bar1$-a.e.\ $\chi\in\X$. If $\bar\tau_*(\chi)\ge1$ for some $\chi$ then $\bar F(\chi)>\|\bar F\|_{C(\X)}$, which is impossible. Hence $\bar\tau_*(\chi)<1$ for $\bar1$-a.e.\ $\chi$. Similarly, $\bar\tau_*(\chi)>0$ for $\bar1$-a.e.\ $\chi$. Therefore, the condition $0\in\partial_\Phi\mathcal{L}$ implies that
	\[
		\bar F(\chi) = \frac{\mu}{2\pi}(2\bar\tau_*(\chi)-1)\sign\phi
		\quad\text{for $\bar1$-a.e.}\quad\chi\in\X.
	\]
	So
	\[
		\bar\tau_*(\chi) = \frac12+\frac\pi\mu\bar F(\chi)\sign\phi
		\qquad\text{and}\qquad
		\bar\Phi(\chi) = \frac\pi\mu |\phi| \bar F(\chi) + \frac12\phi.
	\]
	These equalities holds for a.e.\ $\chi$ wrt $\bar 1$. Since the functions $\bar\Phi(\chi)$ and $\bar\tau_*$ are defined up to a sets of points $\chi$ of $\bar1$-zero measure, we are able to redefine them so that the previous equalities hold for all $\chi$. Therefore, $\bar\Phi(\chi)$ and $\bar\tau_*(\chi)$ become continuous functions on $\X$, i.e. $\Phi,\tau_*\in L_\infty(\SI)$. In particular, we have proved the second part of~\eqref{eq:optimal_phi_L_one_case}.
	
	Now we are able to find $\phi$ from the equation~$s_{D_\infty}(\Phi,\phi)=1$:
	\begin{align*}
		s_{D_\infty}(\Phi,\phi) =& \frac1{2\pi}\int_0^{2\pi}\int_0^1|\Phi(x)-\phi\tau|\D \tau\D x =
			\frac{\sign\phi}{2\pi}\int_0^{2\pi} \left[\int_0^{\tau_*} - \int_{\tau_*}^1\right](\Phi(x)-\phi\tau)\D \tau\D x =\\
		=&\frac{\sign\phi}{2\pi}\int_0^{2\pi} \left[(2\tau_*(x)-1)\Phi(x) - \phi\left(\tau_*^2(x)-\frac12\right)\right]\D x =
			\frac1{2\pi}\int_0^{2\pi}|\phi|\left[ \frac{\pi^2}{\mu^2}F^2(x)+\frac14\right]\D x =|\phi| \left[ \frac{\pi^2}{\mu^2}\langle F^2\rangle + \frac14 \right],		
	\end{align*}
	which proves the first part of~\eqref{eq:optimal_phi_L_one_case} up to the sign of $\phi$.
	
	Let us now use the condition
	\[
		0\in\partial_\phi\mathcal{L}
		\qquad\Leftrightarrow\quad
		-2\pi f \in
		\frac{\mu}{2\pi}\int_0^1\int_0^{2\pi}\tau\sign(\Phi(x)-\phi\tau)\D x\D \tau
	\]
	This integral is easy to compute by Fubini's theorem, since we have shown that $\Phi\in L_\infty(\SI)$:
	\[
		\int_0^1\tau\sign(\Phi(x)-\phi\tau)\D \tau =
		\left(
			\int_0^{\tau_*}\tau\D \tau - \int_{\tau_*}^1\tau\D t
		\right)\sign\phi =
		\frac12(2\tau_*^2-1)\sign\phi =
		\left(
			\frac{\pi^2}{\mu^2}F^2(x) +
			\frac\pi\mu F(x)\sign\phi -
			\frac14
		\right)\sign\phi,
	\]
	Thus
	\[
		-2\pi f =
		\frac{\pi^2}{\mu}\langle F^2\rangle\sign\phi +
		\pi\langle F\rangle -
		\frac\mu4\sign\phi
		\qquad\Leftrightarrow\qquad
		Q(\mu)=\mu^2 - 4\pi(\langle F\rangle + 2f)\mu\sign\phi - 4\pi^2\langle F^2\rangle=0.
	\]
	
Polynomial $Q(\mu)$ has a positive root $\mu_1=\mu_{D_\infty}(F,f)>2\pi\|F\|_{L_\infty(\SI)}$ and another
root $\mu_2$ is non-positive since $4\pi^2\langle F^2\rangle\ge0$. Moreover,  $|\mu_2|<\mu_1$ since $\mu_1|\mu_2|=4\pi^2\langle F^2\rangle\le 4\pi^2\|F\|^2_{L_\infty(\SI)}<\mu_1^2$. Hence the sum of roots is positive, and by Vieta's theorem we obtain
	\[
		\sign\phi = \sign (\langle F\rangle + 2f),
	\]
	which completes the proof of~\eqref{eq:optimal_phi_L_one_case}. Moreover, $Q(\mu)=0$ implies~\eqref{eq:mu_second}.
	
	\medskip
	
	Summarizing, we have proved that if $\mu>2\pi\|F\|_{L_\infty(\SI)}$, then it satisfies~\eqref{eq:mu_second}, and the unique solution $(\Phi,\phi)$ to the auxiliary problem is  given in the theorem statement.
	
	\medskip
	
	It remains to consider the case when the rhs of~\eqref{eq:mu_second} is strictly grater than $2\pi\|F\|_{L_\infty(\SI)}$. In this case, we use~\eqref{eq:mu_second} as a definition for $\mu$ in KKT conditions. Since $\mu>2\pi\|F\|_{L_\infty(\SI)}$, we obtain $\tau_*(x)\in(0;1)$ for a.e.\ $x$ and hence for $(\Phi,\phi)$ given in~\eqref{eq:optimal_phi_L_one_case}, we have $s_{D_\infty}(\Phi,\phi)=1$, $0\in \partial\mathcal{L}(\Phi,\phi)$, and $\langle F,\Phi\rangle + 2\pi\,f\phi=\mu$, i.e.\ the sufficient KKT conditions are satisfied. Hence $\mu=V(F,f)>2\pi\|F\|_{L_\infty(\SI)}$ and so the auxiliary problem has a unique solution as we have already shown.

\end{proof}

\begin{corollary}

	$D_\infty = D_\infty^1 \cap D_\infty^2$ where
	\[
		2\pi D_\infty^1 = \Big\{(F,f):\|F\|_{L_\infty(\SI)}\le 1\Big\}
		\quad\text{and}\quad
		2\pi D_\infty^2 = \Big\{
			(F,f):\big|\langle F\rangle + 2f\big| + \sqrt{(\langle F\rangle + 2f)^2+\langle F^2\rangle}\le 1
		\Big\}.
	\]
	
\end{corollary}

\section{Approximate reach time}\label{sec:reach_time}

\begin{theorem}
\label{thm:T_estimation}
	
	There exists a constant $C>0$ such that for all $(F,f)\in L_\infty(\SI)\times\R$ with $\|F\|_{L_\infty}+|f|^{1/2}\ge C$ the minimal reach time satisfies
	\begin{equation}
	\label{estimate}
		1-C\big(\|F\|_{L_\infty}+|f|^{1/2}\big)^{-1/2}\le
		\frac{T(F,f)}{\max\big\{T_0(F,f),T_1(F,f)\big\}}\le
		1+C\big(\|F\|_{L_\infty}+|f|^{1/2}\big)^{-1}
	\end{equation}
	where $T_0(F,f)=2\pi\|F\|_{L_\infty(\SI)}$ and $T_1(F,f)$ is a unique positive solution to the following ``quadratic'' equation for $T_1$:
	\begin{equation}
	\label{eq:defn_T_1}
		 2\left|
		 	2f+T_1\langle F\rangle
		 \right| =
		 \frac{T_1^2}{2\pi}-2\pi\langle F^2\rangle
	\end{equation}
	
\end{theorem}
\begin{remark}
	
	Thus,
	\begin{equation}
	\label{Bellman_function}
		\T(F,f)\eqdef\max\big\{T_0(F,f),T_1(F,f)\big\}
	\end{equation}
	is an asymptotic Bellman function in our setting. It is interesting that the rhs and the lhs of (\ref{estimate}) have different rates of convergence to 1 as the point $(F,f)$ goes to infinity. In particular, there exists $c_2>0$ such that for any $(F,f)$ with $\|F\|_\infty+|f|^{1/2}>\ge c_2$, one has
	\[
		\frac1{c_2}\big(\|F\|_{L_\infty(\SI)} + |f|^{1/2}\big) \le
		T(F,f) \le
		c_2\big(\|F\|_{L_\infty(\SI)} + |f|^{1/2}\big).
	\]
\end{remark}

\begin{proof}

	Let us take $C>0$ such that $\|(F,f)\|_{L_\infty(\SI)\times\R}\ge c_1$ and $T(F,f)\ge T_1$ where $c_1$ is from Corollary~\ref{cor:reach_time_simple_estimation} and $T_1$ is from Proposition~\ref{prop:D_infty_CTDT_approx} (such a $C$ exists, since $D(T_1)$ is a bounded set).

	Put $T=T(F,f)$ for short. Obviously $(F,f)\not\in D(\tilde T)$ for any $\tilde T<T$, and $(F,f)\in D(T)$ by Corollay~\ref{cor:existance_optimal_control}. We claim that
	\[
		|\mu_{D_\infty}(C(T)(F,f))-1|\le \frac{T_1}T	
	\]
	Indeed, $\mu_{D_\infty}(C(\tilde T)(F,f))$ is a continuous function on $\tilde T$. Thus Corollary~\ref{cor:Minkowski_functional_approx} implies that it is less than $1+T_1/\tilde T$ for $\tilde T\le T$, and $\mu_{D_\infty}(C(T)(F,f))\ge 1-T_1/T$.
	
	Moreover, by using Corollary~\ref{cor:reach_time_simple_estimation}, we obtain
	\begin{equation}
	\label{eq:mu_1_estimate}
		|\mu_{D_\infty}(C(T)(F,f))-1| \le \frac{c_1T_1}{\|F\|_{L_\infty(\SI)} + |f|^{1/2}}.
	\end{equation}
The rest of the  proof goes as follows: we find a solution $T'$ to the equation $\mu_{D_\infty}(C(T')(F,f))=1$, and then we estimate the difference between $T'$ and $T$ using~\eqref{eq:mu_1_estimate}.
	
	Obviously, $\mu_{D_\infty}(C(T')(F,f))$ goes to $0$ as $T'\to\infty$ and it goes to $\infty$ as $T'\to 0$. Since $\mu_{D_\infty}(C(T')(F,f))$ is a continuous function on $T'$, there exists a solution to the equation $\mu_{D_\infty}(C(T')(F,f))=1$. Let us show that this solution is unique. Denote
	\[
		M(b,c,d)=|b+2c|+\sqrt{|b+2c|^2+d}.
	\]
	Obviously, $\mu_{D_\infty}(F,f)=2\pi\max\{\|F\|_\infty,M(\langle F\rangle,f,\langle F^2\rangle)\}$
	
	\begin{lemma}
	\label{lm:monotonicity}
		Fix $(b,c,d)\ne(0,0,0)$. If $d\ge b^2$, then $M(\tau b,\tau^2 c,\tau^2d)$ is strictly increasing in $\tau\ge 0$ and has positive left and right derivatives wrt $\tau$ if $b\ne 2c\tau$ or $d>b^2$.
	\end{lemma}

	\begin{proof}
	
		Define
		\[
			m(\tau)\eqdef\tau\left(
				\big|b + 2\tau c\big| + A(\tau)
			\right).
		\]
		where $A(\tau)=\sqrt{(b + 2\tau c)^2+d}$ for short. Obviously, $\dot m(0)\ge 0$. Suppose, $\tau>0$. Put $\sigma=\sign (b + 2\tau c)$. A direct computation gives
		\[
			\dot m(\tau) =
			\frac{m(\tau)}{\tau}\left(
				1+\frac{2\sigma \tau c}{A(\tau)}
			\right).
		\]
		
		Let us estimate the last factor. Obviously, if $\sign c=\sigma$, then $\dot m(\tau)\ge m(\tau)/\tau$. Suppose that $\sign c=-\sigma$. In this case, $\sign b=\sigma$. We claim that
		\[
			B(\tau)\eqdef d + b^2 + 4bc\tau \ge 0.
		\]
		Indeed, $B(\tau)=d - b^2 + 2b (b + 2c\tau)$. So
		\[
			A(\tau)^2=(b + 2c\tau)^2 + 4c^2\tau^2 = 4c^2\tau^2 + B(\tau)
			\qquad\text{or}\qquad
			A(\tau) = -2\sigma\tau c\sqrt{1+\frac{B(\tau)}{4c^2\tau^2}}.
		\]
		In particular,
		\[
			-2\sigma c\tau\left(1+\frac{B(\tau)}{4c^2\tau^2}\right) =
			A(\tau)\sqrt{1+\frac{B(\tau)}{4c^2\tau^2}}\le
			A(\tau)\left(1+\frac{B(\tau)}{8c^2\tau^2}\right)
			\qquad\text{or}\qquad
			-2\sigma\tau f\le
			A(\tau)\left(1-\frac{B(\tau)}{8c^2\tau^2+2B(\tau)}\right).
		\]
		The last inequality is equivalent to
		\[
			1+\frac{2\sigma c\tau }{A(\tau)} \ge \frac{B(\tau)}{8c^2\tau^2+2B(\tau)}\ge 0.
		\]
		Therefore, $\dot m(\tau)\ge 0$ also in the case $\sign c=-\sigma$. Moreover, the equality $\dot m(\tau)=0$ is possible only if $\sign c=-\sigma$ and $B(\tau)=0$, i.e.\ $d=b^2$ and ($b=0$ or $b+2c\tau=0$). If $b=0$ then $\sign c=\sigma$. Therefore the equality $\dot m(\tau)=0$ may happen only if $b+2c\tau=0$ and $b^2=d$.
	\end{proof}

\begin{remark}
 This lemma is ``qualitatively speaking'' quite  natural, since its geometric analog simply says that reachable sets $D(T)$ are monotone increasing with $T$.
\end{remark}
	
	Lemma~\ref{lm:monotonicity} implies that the equation $\mu_{D_\infty}(\tau F,\tau^2f)=1$ has a unique solution and hence $T'=1/\tau$. The instant $T'$ can be easily found by the exact formula for $\mu_{D_\infty}$ given in Theorem~\ref{thm:main}.
	
	Now, we want to estimate the difference between $T$ and $T'$. Denote
	\[
		K=\Big\{(a,b,c,d):b^2\le d\le a^2,\ a\ge 0,\ 2\pi\max\{a,M(b,c,d)\}=1\Big\}\subset\R^4.
	\]
	Clearly, $K$ is a compact set. Define
	\[
		\beta^+(\tau) = \max_{(a,b,c,d)\in K}\max\big\{\tau a,M(\tau b,\tau^2 c,\tau^2 d)\big\},
		\qquad\text{and}\qquad
		\beta^-(\tau) = \min_{(a,b,c,d)\in K}\max\big\{\tau a,M(\tau b,\tau^2 c,\tau^2 d)\big\}
	\]
	Both functions $\beta^{\pm}$ are strictly monotone by Lemma~\ref{lm:monotonicity}. Both of them are continuous since $M$ is uniformly continuous on any compact. Moreover, $\beta^\pm(0)=0$ and $\beta^\pm(\tau)\to\infty$ as $\tau\to\infty$.
	
	We need to study behavior of $\beta^\pm$ near $\tau=1$. By Lemma~\ref{lm:monotonicity}, function $M(\tau b,\tau^2 c,\tau^2 d)$ may have zero derivative at $\tau=1$ only if $b=-2c$ and $d=b^2$. Hence for $(a,b,c,d)\in K$, the left or right derivative of the function $\max\{\tau a,M(\tau b,\tau^2 c,\tau^2d)\}$ may vanish only when $b=-2c$, $d=b^2$, $2\pi a\le 1$, $d\le a^2$, and $2\pi M(b,c,d)=1$. These equations have only 2 possible solutions: $a_0=\frac{1}{2\pi}$, $b_0^\pm=\pm a_0$, $c_0^\pm=-\frac12 b_0^{\pm}$, and $d_0=a_0^2$.
	
	We start with investigation of $\beta^-(\tau)$ behavior as $\tau\to 1+0$. Outside a neighborhood of these two points $(a_0,b_0^\pm,c_0^\pm,d_0)$, right derivatives of the function $M(\tau b,\tau^2 c,\tau^2d)$ wrt $\tau$ at $\tau=1$ is positive and separated from zero. Obviously, derivative of $\tau a$ is positive and separated from 0 in a neighborhood of points in $K$ with $a=0$. Moreover we may chose these neighborhoods in such a way that they do not intersect each other. Hence $2\pi\max\{\tau a,M(\tau b,\tau^2 c,\tau^2d)\}\ge 1+c(\tau-1)$ for some constant $c$ and $\tau\ge 1$ that is close to $1$. Therefore $\beta^-(\tau)\ge\frac1{2\pi}+c(\tau-1)$ for the same $\tau$.
	
	Behavior of $\beta^+(\tau)$ for $\tau\to 1-0$ is different. Outside a neighborhood of points $(a_0,b_0^\pm,c_0^\pm,d_0)$ the left derivative at $\tau=1$ is also positive and separated from 0, but now $\max\{\tau a_0,M(\tau b_0^\pm,\tau^2 c_0^\pm,\tau^2d_0)\}=M(\tau b_0^\pm,\tau^2 c_0^\pm,\tau^2d_0)$ for $\tau\to 1-0$. So $\beta^+(\tau)$ does not have a linear growth as $\tau\to 1-0$. We claim that $\beta^+(\tau)\ge \frac1{2\pi}-c(1-\tau)^2$ for some constant $c>0$ and $\tau<1$ that is close to 1. Indeed, left second derivative of $M(\tau b_\pm,\tau^2 c_\pm,\tau^2d_\pm)$ wrt $\tau$ at $\tau=1$ is negative.
	
	Consider two equations $\beta^\pm(\tau)=\frac1{2\pi}+\varepsilon$. Each of them has a unique solution $\tau^\pm(\varepsilon)$ for any $\varepsilon\ge -1$. Both $\tau^\pm(\varepsilon)$ are continuous strictly increasing functions, and $\tau^\pm(0)=1$. Moreover, $\tau^-(\varepsilon)\le1+c\varepsilon$ as $\varepsilon\to +0$ and $\tau^+(\varepsilon)=1-c\varepsilon^{1/2}$ as $\varepsilon\to -0$.
	
	We know that
	\[
		(a',b',c',d')\eqdef
		(T'^{-1}\|F\|_{L_\infty(\SI)},T'^{-1}\langle F\rangle,T'^{-2}f,T'^{-1}\langle F^2\rangle)
		\in K
	\]
	by definition of $T'$. Put $\tau=T'/T$. Hence, estimate~\eqref{eq:mu_1_estimate} implies
	\[
		\frac1{2\pi}-\varepsilon\le
		\max\big\{\tau a',M(\tau b',\tau^2c',\tau^2d')\big\}\le
		\frac1{2\pi}+\varepsilon
		\qquad\text{where}\qquad
		\varepsilon=\frac{c_1T_1}{\|F\|_{L_\infty(\SI)} + |f|^{1/2}}.
	\]
	Therefore, $\beta^-(\tau)\le \frac1{2\pi}+\varepsilon$, and $\beta^+(\tau)\ge \frac1{2\pi}-\varepsilon$. We conclude that $\tau^+(-\varepsilon)\le\tau\le \tau^-(\varepsilon)$.

\end{proof}

\begin{remark}
\label{rm:bad_curve}
	
	Note that the lhs estimate in~\eqref{estimate} deteriorates only near points $(F,f)$ that satisfy $\|F\|_{L_\infty(\SI)}=\T/2\pi$, $\langle F\rangle=\pm\T/2\pi$, $f=\mp \T^2/4\pi$, and $\langle F^2\rangle=\T^2/4\pi^2$. These conditions obviously imply that $F$ identically equal to some constant, which then must be equal to $\pm\T/2\pi$. So points $F=\pm T/2\pi$ and $f=\mp T^2/4\pi$ for $T\in\R_+$ form two 1-dimensional curves in the state space $L_\infty(\SI)\times\R$. If we decide to work outside these curves, e.g.\ chose their neighborhood that is invariant under scaling operators $C(T)$ and consider only pairs $(F,f)$ that does not belong to the neighborhood, then the estimate~\eqref{estimate} becomes much better:
	\[
		T(F,f)=\T(F,f)(1+O(\|F\|_{L_\infty}+|f|^{1/2})^{-1})
	\]
	
\end{remark}
	
\section{Analysis of asymptotically optimal control}
\label{sec:synthesis}

By definition, the asymptotically optimal control in our setting is a control which allows to come from zero to a given state $(F,f)$ in time $\T(F,f)$, which is asymptotically equivalent to the minimum time $T(F,f)$: $$\frac{\T(F,f)}{T(F,f)}\mbox{ tends to 1 as }(F,f)\mbox{ tends to }\infty. $$

We turn to the governing equations \eqref{eq:basic_system2} and convert them into something as explicit as possible. For simplicity assume that the motion starts at time $t=0$. Then it follows from  \eqref{eq:basic_system2} that if $s\in[0,2\pi]$
\begin{equation}
\label{eq:F_motion}
	F(x)|_{t=s}=F(x-s)|_{t=0}+
	\begin{cases}
		u(2\pi-x),&\text{if}\quad0<x<s;\\
		0,&\text{otherwise.}
	\end{cases}
\end{equation}
In particular, for $t=2\pi$ we obtain
\begin{equation}
\label{eq:F_period}
    F(x)|_{t=2\pi}=F(x)|_{t=0}+ u(2\pi-x),
\end{equation}
Behavior of $f$ in time is more complicated, but we do need its behavior only for the case when $u\equiv 1$ or $u\equiv -1$ during the period. In this case, it is easy to compute
\begin{equation}
\label{eq:f_period}
	f|_{t=2\pi} = f|_{t=0} -2\pi\left(\langle F\rangle +\frac12\langle u\rangle\right).
\end{equation}

We are going to study how the transformation $(F,f)|_{t=0}\mapsto(F,f)|_{t=2\pi}$ affects the asymptotic Bellman function $\T=\max(T_0,T_1)$ given in~\eqref{Bellman_function} for the case when $\|F\|_\infty+|f|^{1/2}$ is sufficiently large. This, of course, depends on the choice of control $u(t),\,\, t\in[0,2\pi]$.

First suppose that $\T=T_0$. Following~\cite{Ovseevich2018}, we put $u(t)=\sign F(2\pi-t)|_{0}$. Equation~\eqref{eq:F_motion} implies
\[
	F(x)|_t=F(x-t)|_{t=0}+\sign F(x-t)|_{t=0},
\]
Hence
\[
	F(x)|_{t=2\pi}=F(x)|_0+ \sign F(x)|_{t=0}.
\]
In particular,
\[
  \big|F(x)|_{t=2\pi}\big|=\big|F(x)|_{t=0}\big|+1
\]
provided that $\big|F(x)|_{0}\big|\ne 0$. In particular, the constituent $T_0(F)=2\pi\|F\|_{L_\infty(\SI)}$ of the asymptotic Bellman function increases by $2\pi$ during the period $[0,2\pi]$. Therefore, the asymptotic Bellman function increases by $2\pi$ at least. This is certainly the best behavior of $\T$ we can expect.

Now we turn to the case $\T=T_1$, i.e. $T_1\ge 2\pi\|F\|_{L_\infty(\T)}$, where we write $F$ instead of $F|_{t=0}$  for short. Note that in the considering case $\T-2\pi|\langle F\rangle|\ge 0$, but it may still happen that the equality if fulfilled. This situation may only appear if $\T=2\pi\|F\|_{L_\infty(\SI)}$ and $F\equiv\mathrm{const}=\pm\T/2\pi$. In this rare case, since $\langle F\rangle^2=\langle F^2\rangle$, \eqref{eq:defn_T_1} implies $f=\mp \T^2/2\pi$. So the described situation may happen only on two 1-dimensional curves in the state space $L_{\infty}(\SI)\times\R$, which identically coincide with the curves described in Remark~\ref{rm:bad_curve}.

Lat us proceed to developing asymptotically optimal control. As it is pointed out in Proposition~\ref{prop:DT_support_func}, we should consider $\sign(\Phi(T-t)+(T-t)\varphi)$ as a control. This control seems a little bit complicated. Fortunately, when constructing asymptotically optimal control, we are not obligated to use exact formulae, since small errors cannot destroy asymptotic optimality. Hence, for simplicity, we assume that $u(t)$ is $1$ or $-1$ during the whole period\footnote{If $2f+\T\langle F\rangle=0$ then both choices $\sigma=1$ and $\sigma=-1$ are allowed, but the control $u(t)$ must be constant during the period.}:
\[
u(t)\equiv\sigma\eqdef\sign\phi(C(\T)(F,f)|_{t=0})=
\sign\left(\T^{-1}\langle F\rangle+2\T^{-2}f\right)=
\sign(2f+\T\langle F\rangle)
\]
where we again write $F$ and $f$ instead of $F|_{t=0}$ and $f|_{t=0}$ for short and the previous expression for the sing of $\phi$ is taken from Theorem~\ref{thm:main}.

Denote by $F'=F_{t=2\pi}$ and $f'=f|_{t=2\pi}$. Then
\[
	\begin{array}{lcl}
		2|2f+\T\langle F\rangle|=\frac{\T^2}{2\pi}-2\pi\langle F^2\rangle&\text{ at }&t=0,\\[.5em]
		2|2f'+T'_1\langle F'\rangle|=\frac{{T'_1}^2}{2\pi}-2\pi\langle F'^2\rangle&\text{ at }&t=2\pi,
	\end{array}
\]
where we denote $T'_1=T_1(F',f')$ for short. Obviously $\T(F',f')\ge T'_1$.

Now, we verify that $T'_1=\T+2\pi$, which trivially implies that the chosen control $u$ is an asymptotically optimal. So according to~\eqref{eq:F_period} and~\eqref{eq:f_period}, we have
\begin{equation}
\label{eq:Ff_after_period}
	\begin{gathered}
		\langle F'\rangle =\langle F\rangle+\sigma\\
		f' =f-2\pi\big(\langle F\rangle+ {\textstyle {\frac12}}\sigma\big)
	\end{gathered}
\end{equation}
Hence
\[
  \langle F'^2\rangle =  \langle F^2\rangle +2\sigma\langle F\rangle+1.
\]

On the one hand, since $\T-2\pi\sigma\langle F\rangle\ge0$, we have

\begin{equation}
\label{F22}
	2\big|2f'+(\T+2\pi)\langle F'\rangle\big|=
	2\big|(2f+\T\langle F\rangle)+\sigma(\T-2\pi\sigma\langle F\rangle)\big|=
	2\big|2f+\T\langle F\rangle\big|+2(\T-2\pi\sigma\langle F\rangle).
\end{equation}
On the other hand,
\begin{equation}
\label{F222}
	\frac{(\T+2\pi)^2}{2\pi}-2\pi\langle F'^2\rangle=
	\frac{\T^2}{2\pi}-2\pi\langle F^2\rangle+2(\T-2\pi\sigma\langle F\rangle).
\end{equation}
Note, that both~\eqref{F22} and~\eqref{F222} contain the identical term $2(\T-2\pi\sigma\langle F\rangle)$ in the rhs. Now we recall that $\T=T_1$, which shows by~\eqref{eq:defn_T_1} that
\[
	2\left|2f'+(\T+2\pi)\langle F'\rangle\right|=
	\frac{(\T+2\pi)^2}{2\pi}-2\pi\langle F'^2\rangle.
\]
Therefore $\T+2\pi$ and $T'_1$ are solutions to the same equation, and hence they must coincide due to Lemma~\ref{lm:monotonicity}. Therefore
\[
	\T(F',f') \ge T'_1 = \T(F,f)+2\pi.
\]
This is again an asymptotically optimal behavior.

\medskip

We also decided to add a brief explanation of the case $2f+\T\langle F\rangle=0$. For any choice of $\sigma=\pm1$, after the $2\pi$-period we have $2f'+2\pi\langle F'\rangle = 2f+\T\langle F\rangle + \sigma(\T-2\pi\sigma\langle F\rangle)$. So if $T>2\pi|\langle F\rangle|$, then after a $2\pi$-period we will have $2f'+2\pi\langle F'\rangle\ne 0$ and $\sigma=\sign (2f'+2\pi\langle F'\rangle)$. But in the case of any of two described in Remark~\ref{rm:bad_curve} curves, if $\sigma=\sign\langle F\rangle$, then after the $2\pi$-period we will stay on the curve and have $2f'+2\pi\langle F'\rangle=0$. The opposite choice $\sigma=-\sign\langle F\rangle$ immediately throw the point out of the curve and we again have $2f'+2\pi\langle F'\rangle\ne 0$ and $\sigma=\sign(2f'+2\pi\langle F'\rangle)$. This situation is completely similar to the very well known optimal synthesis in the Pontryagin time minimization problem: $T\to\min$, $|\ddot x|\le 1$, $x\in\R$.

\medskip

Thanks to Theorem~\ref{thm:T_estimation} we know that $T=\T(1+O(\|F\|_\infty+|f|^{1/2})^{-1/2})$. The suggested asymptotically optimal control can double function $\T$ in time $\T$. Therefore, it will double $T$ in time $T(1+O(\|F\|_\infty+|f|^{1/2})^{-1/2})$. The absolutely (but not just asymptotically) optimal result is doubling in time $T$. So the designed asymptotically optimal control works pretty well if $\|F\|_\infty+|f|^{1/2}$ is large enough.

Moreover, if we apply the designed control in backward direction, it allows to halve $T$ in time $\frac12T(1+O(\|F\|_\infty+|f|^{1/2})^{-1/2})$. Repeating this $\sim \log_2T$ times one can reach a fixed bounded neighborhood of $(0,0)$ in time $T(1+O(\|F\|_\infty+|f|^{1/2})^{-1/2})$, which gives an asymptotically optimal solution to the finite-time stabilization problem:
\[
	\begin{gathered}
	T\to\min\qquad\\
	\begin{aligned}
		F|_{t=0}=F_0;\quad&\quad f|_{t=0}=F_0;\\
		F|_{t=T}=0;\quad&\quad f|_{t=T}=0.
	\end{aligned}
	\end{gathered}
\]

\section{Funding}

The work of L.V.~Lokutsievskiy is supported by the Russian Science Foundation under grant 20-11-20169 and performed in Steklov Mathematical Institute of the Russian Academy of Sciences. The work of A.I. Ovseevich was supported by the Russian Science Foundation under grant 21-11-00151 and performed in Ishlinsky Institute for Problems in Mechanics of the Russian Academy of Sciences. Sections 1,3,5,7, and 9 were written by L.V.~Lokutsievskiy. Sections 2,4,6,8, and 10 were written by A.I.~Ovseevich. All results in this paper are products of authors collaborative work.

\end{document}